\documentclass[a4paper,12pt]{article}
\usepackage[latin1]{inputenc}
\usepackage{amssymb}

\usepackage{times}
\newtheorem{The}{Theorem}[section]
\newtheorem{Pro}[The]{Proposition}
\newtheorem{Deff}[The]{Definition}
\newtheorem{Lem}[The]{Lemma}

\newtheorem{Cor}[The]{Corollary}

\title{Decision Problems for Recognizable Languages of Infinite Pictures} 

\author{Olivier Finkel   \\{\it   Equipe de Logique Math\'ematique}
  \\ CNRS et  Universit\'e Paris Diderot Paris 7
 \\ UFR de Math\'ematiques case 7012, site Chevaleret,\\75205 Paris Cedex 13, 
 France.\\ finkel@logique.jussieu.fr}

\date{}
\begin{document}

\newcommand{\om}{\omega}
\newcommand{\Si}{\Sigma}
\newcommand{\Sis}{\Sigma^\star}
\newcommand{\Sio}{\Sigma^\omega}
\newcommand{\nl}{\newline}
\newcommand{\lra}{\leftrightarrow}
\newcommand{\fa}{\forall}
\newcommand{\ra}{\rightarrow}
\newcommand{\orl}{ $\omega$-regular language}

\newcommand{\Ga}{\Gamma}
\newcommand{\Gas}{\Gamma^\star}
\newcommand{\Gao}{\Gamma^\omega}
\newcommand{\ite}{\item}

\newcommand{\ol}{$\omega$-language}

\newcommand{\tla}{\twoheadleftarrow}
\newcommand{\de}{deterministic }
\newcommand{\proof}{\noindent          {\bf Proof.} }
\newcommand {\ep}{\hfill $\square$}

\maketitle

\begin{abstract}
\noindent   Altenbernd, Thomas and 
W\"ohrle have considered in \cite{ATW02} acceptance of 
languages of infinite two-dimensional words (infinite pictures) by finite tiling systems,  
with the usual acceptance conditions, such as  the 
B\"uchi and Muller ones,  firstly  used for infinite words.  
Many classical decision problems are studied in formal language theory and in automata theory and arise now naturally about 
recognizable languages of infinite pictures. 
We first review in  this paper some recent results of \cite{Fink-tilings} where 
we gave  the exact degree of numerous   undecidable problems for B\"uchi-recognizable 
languages of infinite pictures, which are actually  located  at the first or at the 
second level of the analytical hierarchy, and ``highly undecidable". 
Then we prove here some more (high) undecidability results. We first show that it is  $\Pi_2^1$-complete to determine whether a given 
  B\"uchi-recognizable  languages of  infinite pictures  is unambiguous. Then we investigate cardinality problems. 
Using recent results of  \cite{FL-TM}, we prove that it is 
$D_2(\Sigma_1^1)$-complete to determine whether a given B\"uchi-recognizable language of infinite pictures is countably infinite, and that 
 it is $\Sigma_1^1$-complete to determine whether a given B\"uchi-recognizable language of infinite pictures is uncountable. 
Next we consider complements of recognizable languages of infinite pictures.   Using some results of Set Theory, 
we show that the cardinality of the complement of a B\"uchi-recognizable language of infinite pictures 
may depend on the  model of the axiomatic system {\bf ZFC}.  We prove that the problem to determine 
whether the complement of a given B\"uchi-recognizable language of infinite pictures is countable (respectively, uncountable) 
 is in the class $\Si_3^1 \setminus (\Pi_2^1 \cup \Si_2^1)$ (respectively, in the class $\Pi_3^1 \setminus (\Pi_2^1 \cup \Si_2^1)$). 
\end{abstract}


\section{Introduction}

Languages of infinite words accepted by finite automata were first studied by B\"uchi 
to prove the decidability of the monadic second order theory of one successor
over the integers.  Since then regular $\om$-languages have been much studied and  many applications have been  found  for specification and verification 
of non-terminating systems, 
see \cite{Thomas90,PerrinPin} for many results and references. 

\bigskip \noindent                   
Altenbernd, Thomas and 
W\"ohrle have considered in \cite{ATW02} acceptance of 
languages of infinite two-dimensional words (infinite pictures) by finite tiling systems,  
with the usual acceptance conditions, such as  the 
B\"uchi and Muller ones,  firstly  used for infinite words.
This way they extended both  the classical theory of  $\om$-regular languages and the classical theory of recognizable languages of finite pictures, 
  \cite{Giammarresi-Restivo},  to the case of infinite pictures.  

\bigskip \noindent                    Many classical decision problems are studied in formal language theory and in automata theory and arise now naturally about 
recognizable languages of infinite pictures. 

\bigskip \noindent                    In a recent paper,  we gave  the exact degree of numerous   undecidable problems for B\"uchi-recognizable 
languages of infinite pictures.   In particular, the non-emptiness and the infiniteness problems are $\Si_1^1$-complete, and the  universality 
problem, the inclusion problem, the equivalence problem,  the complementability problem,  and the determinizability problem,  are all $\Pi_2^1$-complete.  
These decision problems are then located  at the first or at the 
second level of the analytical hierarchy, and ``highly undecidable". 
This gave new natural examples of decision problems located at the first or at the second level of the analytical hierarchy.  

\bigskip \noindent                   Here we first review some of these results, and we study new decision problems, obtaining new results of high undecidability. 

\bigskip \noindent                   We first consider the notion of unambiguous B\"uchi tiling system, and of unambiguous B\"uchi-recognizable language of infinite pictures. We 
show that every language of infinite pictures which is accepted by an unambiguous B\"uchi tiling system is a Borel set. As a corollary this shows the existence 
of inherently ambiguous B\"uchi-recognizable language of infinite pictures. 
Then we use this result to prove that it is $\Pi_2^1$-complete to determine whether a given B\"uchi-recognizable language of infinite pictures is 
unambiguous. 

\bigskip \noindent                   Next we study cardinality problems. Using recent results of Finkel and Lecomte in \cite{FL-TM}, we first show that it is 
$D_2(\Sigma_1^1)$-complete to determine whether a given B\"uchi-recognizable language of infinite pictures is countably infinite, where 
$D_2(\Sigma_1^1)$ is the class of $2$-differences of $\Sigma_1^1$-sets, i.e. the class 
 of sets which are intersections of a $\Sigma_1^1$-set and of a $\Pi_1^1$-set.  
And it is $\Sigma_1^1$-complete to determine whether a given B\"uchi-recognizable language of infinite pictures is uncountable. 

\bigskip \noindent                   Then we consider the complements of  B\"uchi-recognizable languages of infinite pictures. 
By using some results of Set Theory, we show that the cardinality of the complement of a B\"uchi-recognizable language of infinite pictures 
may depend on the actual model of the axiomatic system {\bf ZFC}.  
We prove that one can effectively construct a B\"uchi tiling system ${\cal T}$ accepting a language $L \subseteq \Si^{\om, \om}$, whose complement 
is  $L^-=\Si^{\om, \om}-L$, such that: 
\begin{enumerate}
\item There is a model $V_1$ of  {\bf ZFC} in which      $L^-$ is countable. 
\item There is a model $V_2$ of  {\bf ZFC} in which     $L^-$ has cardinal $2^{\aleph_0}$. 
\item There is a model $V_3$ of  {\bf ZFC} in which      $L^-$ has cardinal $\aleph_1$ with 
$\aleph_0<\aleph_1<2^{\aleph_0}$.
\end{enumerate}

\noindent          Then,  using the proof of this result and Schoenfield's Absoluteness Theorem, we prove that the problem to determine 
whether the complement of a given B\"uchi-recognizable language of infinite pictures is countable (respectively, uncountable) 
 is in the class $\Si_3^1 \setminus (\Pi_2^1 \cup \Si_2^1)$ (respectively, in the class $\Pi_3^1 \setminus (\Pi_2^1 \cup \Si_2^1)$). 
This shows that natural cardinality problems are  actually 
located at the {\bf third level} of the analytical hierarchy.

\bigskip \noindent                   The paper is organized as follows. We recall in Section $2$ the notions of tiling systems and of recognizable languages of pictures.
In section $3$, we recall the definition of the analytical hierarchy on subsets of $\mathbb{N}$. The definitions of the Borel hierarchy and of analytical sets 
of a Cantor space, along with their effective counterparts, are given in Section $4$. Some notions of Set Theory, which are useful in the sequel, 
are exposed in Section $5$. 
We study decision problems in Section $6$, proving new results. Some concluding remarks are  given in Section $7$.

\section{Tiling Systems}

\noindent          We assume   the reader to be familiar with the theory of formal ($\om$)-languages  
\cite{Thomas90,Staiger97}.
We recall usual notations of formal language theory. 
\nl  When $\Si$ is a finite alphabet, a {\it non-empty finite word} over $\Si$ is any 
sequence $x=a_1\ldots a_k$, where $a_i\in\Sigma$ 
for $i=1,\ldots ,k$ , and  $k$ is an integer $\geq 1$. The {\it length}
 of $x$ is $k$, denoted by $|x|$.
 The {\it empty word} has no letter and is denoted by $\lambda$; its length is $0$. 
 $\Sis$  is the {\it set of finite words} (including the empty word) over $\Sigma$.
 \nl  The {\it first infinite ordinal} is $\om$.
 An $\om$-{\it word} over $\Si$ is an $\om$ -sequence $a_1 \ldots a_n \ldots$, where for all 
integers $ i\geq 1$, ~
$a_i \in\Sigma$.  When $\sigma$ is an $\om$-word over $\Si$, we write
 $\sigma =\sigma(1)\sigma(2)\ldots \sigma(n) \ldots $,  where for all $i$,~ $\sigma(i)\in \Si$,
and $\sigma[n]=\sigma(1)\sigma(2)\ldots \sigma(n)$  for all $n\geq 1$ and $\sigma[0]=\lambda$.
\nl 
 The usual concatenation  of two finite words $u$ and $v$ is 
denoted $u.v$ (and sometimes just $uv$). This product is extended to the product of a 
finite word $u$ and an $\om$-word $v$: the infinite word $u.v$ is then the $\om$-word such that:
\nl $(u.v)(k)=u(k)$  if $k\leq |u|$ , and 
 $(u.v)(k)=v(k-|u|)$  if $k>|u|$.
\nl   
 The {\it set of } $\om$-{\it words } over  the alphabet $\Si$ is denoted by $\Si^\om$.
An  $\om$-{\it language} over an alphabet $\Sigma$ is a subset of  $\Si^\om$.  

\bigskip \noindent                   We now  define two-dimensional words, i.e. pictures.
\nl  Let $\Si$ be a finite alphabet, let $\#$ be a letter not in $\Si$ and let 
$\hat{\Si}=\Si \cup \{\#\}$. If $m$ and $n$ are two positive integers  or if $m=n=0$,  
a   picture of size $(m, n)$ over $\Si$ 
is a function $p$ from $\{0, 1, \ldots , m+1\} \times \{0, 1, \ldots , n+1\}$ 
into $\hat{\Si}$ such that $p(i, j)=\#$ if  $i\in \{0, m+1\}$ or $j\in \{0, n+1\}$ and 
$p(i, j) \in \Si$ otherwise. The empty picture 
is the only picture of size $(0, 0)$ and is denoted by $\lambda$. Pictures of  size 
$(n, 0)$ or $(0, n)$, for $n>0$, are not defined. $\Si^{\star, \star}$ is the set of 
pictures over $\Si$. A picture language $L$ is a subset of $\Si^{\star, \star}$. 
The research on  picture languages was firstly
motivated by the problems arising in pattern recognition and image
processing, a survey on the  theory of picture languages may be found
in  \cite{Giammarresi-Restivo}. 

\bigskip \noindent                   An $\om$-picture over $\Si$ 
is a function $p$ from $\om \times \om$ into $\hat{\Si}$ such that $p(i, 0)=p(0, i)=\#$ 
for all $i\geq 0$ and $p(i, j) \in \Si$ for $i, j >0$. $\Si^{\om, \om}$ is the set of 
$\om$-pictures over $\Si$. An $\om$-picture language $L$ is a subset of $\Si^{\om, \om}$. 
\nl   For $\Si$ a finite alphabet we call  $\Si^{\om^2}$  the set of functions 
from $\om \times \om$ into $\Si$. So the set $\Si^{\om, \om}$ of $\om$-pictures over 
$\Si$ is a strict subset of $\hat{\Si}^{\om^2}$.  

\bigskip \noindent                   We shall say that, for each integer $j\geq 1$,  the $j^{th}$ row of an $\om$-picture 
$p\in \Si^{\om, \om}$ is the infinite word $p(1, j).p(2, j).p(3, j) \ldots$ over $\Si$ 
and the $j^{th}$ column of $p$ is the infinite word $p(j, 1).p(j, 2).p(j, 3) \ldots$ 
over $\Si$. 
\nl As usual,  one can imagine that,  for integers $j > k \geq 1$, the $j^{th}$ column of $p$ 
is on the right of the $k^{th}$ column of $p$ and that the $j^{th}$ row of $p$ is 
``above" the $k^{th}$ row of $p$. 

\bigskip \noindent                   We introduce now (non deterministic)  tiling systems as in the paper \cite{ATW02}. 
\nl A tiling system is a tuple ${\cal A}$=$(Q, \Si, \Delta)$, where $Q$ is a finite set 
of states, $\Si$ is a finite alphabet, $\Delta \subseteq (\hat{\Si} \times Q)^4$ is a finite set 
of tiles. 
\nl A B\"uchi tiling system is a pair $({\cal A},$$ F)$ 
 where ${\cal A}$=$(Q, \Si, \Delta)$ 
is a tiling system and $F\subseteq Q$ is the set of accepting states. 
\nl A Muller tiling system is a pair $({\cal A}, {\cal F})$ 
where ${\cal A}$=$(Q, \Si, \Delta)$ 
is a tiling system and ${\cal F}$$\subseteq 2^Q$ is the set of accepting sets of states.

\bigskip \noindent                    Tiles are denoted by 
$
\left ( \begin{array}{cc}  (a_3, q_3) & (a_4, q_4)
\\ (a_1, q_1) & (a_2, q_2) \end{array} \right ) \mbox{ with } a_i \in \hat{\Si} \mbox{ and } 
q_i \in Q, 
$

\bigskip \noindent                   and in general, over an alphabet $\Ga$, by 
$
\left ( \begin{array}{cc} b_3 & b_4 
\\ b_1 & b_2 \end{array} \right ) \mbox{  ~~~~~~with  }  b_i \in \Ga . 
$

\noindent          A combination of tiles is defined  by: 
\begin{displaymath}
\left ( \begin{array}{cc}  b_3 & b_4
\\ b_1 & b_2 \end{array} \right ) \circ  
\left ( \begin{array}{cc} b'_3 & b'_4
\\ b'_1 & b'_2  \end{array} \right ) = 
\left ( \begin{array}{cc} (b_3, b'_3) & (b_4, b'_4) 
\\ (b_1, b'_1) & (b_2, b'_2) \end{array} \right ) 
\end{displaymath}

\noindent          A run of a tiling system ${\cal A}$=$(Q, \Si, \Delta)$ over a (finite) 
picture $p$ of size $(m, n)$ over $\Si$ 
is a mapping $\rho$  from $\{0, 1, \ldots , m+1\} \times \{0, 1, \ldots , n+1\}$ 
into $Q$ such that for all $(i, j) \in \{0, 1, \ldots , m\} \times \{0, 1, \ldots , n\}$ 
with $p(i, j)=a_{i, j}$ and $\rho(i, j)=q_{i, j}$ we have 
\begin{displaymath}
\left ( \begin{array}{cc}  a_{i, j+1} & a_{i+1, j+1} 
\\  a_{i, j}  & a_{i+1, j}   \end{array} \right ) \circ  
\left ( \begin{array}{cc} q_{i, j+1} & q_{i+1, j+1} 
\\ q_{i, j}  & q_{i+1, j}   \end{array} \right ) \in \Delta . 
\end{displaymath}

\noindent          A run of a tiling system ${\cal A}$=$(Q, \Si, \Delta)$ over an 
$\om$-picture $p \in \Si^{\om, \om}$ 
is a mapping $\rho$  from $\om \times \om$ 
into $Q$ such that for all $(i, j) \in \om \times \om$ 
with $p(i, j)=a_{i, j}$ and $\rho(i, j)=q_{i, j}$ we have 
\begin{displaymath}
\left ( \begin{array}{cc}  a_{i, j+1} & a_{i+1, j+1} 
\\  a_{i, j}  & a_{i+1, j}   \end{array} \right ) \circ  
\left ( \begin{array}{cc} q_{i, j+1} & q_{i+1, j+1} 
\\ q_{i, j}  & q_{i+1, j}   \end{array} \right ) \in \Delta . 
\end{displaymath}

\noindent          We  now recall  acceptance of finite or  infinite pictures by tiling systems:

\begin{Deff} Let ${\cal A}$=$(Q, \Si, \Delta)$ 
be a tiling system, $F\subseteq Q$ and ${\cal F}$$\subseteq 2^Q$. 
\begin{itemize}

\ite
 The picture language recognized by ${\cal A}$ 
is the set of pictures $p \in \Si^{\star, \star}$ such that there is some run $\rho$ of 
${\cal A}$ on $p$. 

\ite The $\om$-picture language  B\"uchi-recognized  
by 
$({\cal A},$$ F)$ 
is the set of $\om$-pictures $p \in \Si^{\om, \om}$ such that there is some run $\rho$ of 
${\cal A}$ on $p$ and $\rho(v) \in F$   for infinitely many $v\in \om^2$.    It is denoted by  
$L^B(({\cal A},$$ F))$.

\ite The $\om$-picture language Muller-recognized  by $({\cal A}, {\cal F})$ is 
the set of $\om$-pictures $p \in \Si^{\om, \om}$ such that there is some run $\rho$ of 
${\cal A}$ on $p$ and $Inf(\rho) \in {\cal F}$ where $Inf(\rho)$ is the set of states 
occurring infinitely often in $\rho$. It is denoted by $L^M(({\cal A},$$ {\cal F}))$.
\end{itemize}
\end{Deff}

\noindent          Notice that an $\om$-picture language $L \subseteq \Si^{\om, \om}$ is 
recognized by a B\"uchi tiling system if and only if it 
is recognized by a  Muller tiling system, \cite{ATW02}.  
\nl We shall denote $TS(\Si^{\om, \om})$ the class of languages $L \subseteq \Si^{\om, \om}$ 
which are recognized by some  B\"uchi (or Muller) tiling system.

\section{Recall of Known Basic Notions}

\subsection{The Analytical Hierarchy}

\noindent          
 The set of natural numbers is denoted by $\mathbb{N}$ and the set of all mappings  from $\mathbb{N}$ into $\mathbb{N}$ will 
be denoted by ${\cal F}$. 

\bigskip \noindent                   
We assume the reader to be familiar with the arithmetical   hierarchy on subsets of  $\mathbb{N}$.   We now recall   the notions
 of analytical hierarchy and of complete sets for classes of this hierarchy which may be found in \cite{rog}.

\begin{Deff}
Let  $k, l >0$ be some  integers. $\Phi$ is a partial recursive function of $k$ function variables and $l$ number variables if there exists $z\in \mathbb{N}$ 
such that for any $(f_1, \ldots , f_k, x_1, \ldots , x_l) \in {\cal F}^k \times \mathbb{N}^l$, we have 
$$\Phi (f_1, \ldots , f_k, x_1, \ldots , x_l) = \tau_z^{f_1, \ldots , f_k}(x_1, \ldots , x_l),$$
\noindent          where the right hand side is the output of the Turing machine with index $z$ and oracles $f_1, \ldots , f_k$ over the input $(x_1, \ldots , x_l)$. 
For $k>0$ and $l=0$, $\Phi$ is a partial recursive function if, for some $z$, 
$$\Phi (f_1, \ldots , f_k) = \tau_z^{f_1, \ldots , f_k}(0).$$
\noindent          The value $z$ is called the Gödel number or index for $\Phi$. 
\end{Deff}

\begin{Deff}
Let $k, l >0$ be some  integers and $R  \subseteq {\cal F}^k \times \mathbb{N}^l$. The relation $R$ is said to be a recursive relation of $k$ 
function variables and $l$ number variables if its characteristic function is recursive. 
\end{Deff}

\noindent          We now define analytical subsets of $\mathbb{N}^l$.

\begin{Deff}
A subset $R$ of $\mathbb{N}^l$ is analytical if it is recursive or if there exists a recursive set $S  \subseteq {\cal F}^m \times \mathbb{N}^n$, with 
$m\geq 0$ and $n\geq l$, such that 
$$R = \{ (x_1, \ldots , x_l) \mid (Q_1s_1)(Q_2s_2) \ldots (Q_{m+n-l}s_{m+n-l}) S(f_1, \ldots , f_m, x_1, \ldots , x_n) \}, $$
\noindent          where $Q_i$ is either $\fa$ or $\exists$ for $1 \leq i \leq m+n-l$, and where $s_1, \ldots , s_{m+n-l}$ are $f_1, \ldots , f_m, x_{l+1}, \ldots , x_n$ in 
some order. 
\nl The expression $(Q_1s_1)(Q_2s_2) \ldots (Q_{m+n-l}s_{m+n-l}) S(f_1, \ldots , f_m, x_1, \ldots , x_n)$ is called a predicate form for $R$. A
quantifier applying over a function variable is of type $1$, otherwise it is of type $0$. 
In a predicate form the (possibly empty) sequence of quantifiers, indexed by their type, is called the prefix of the form. The reduced prefix is the sequence of 
quantifiers obtained by suppressing the quantifiers of type $0$ from the prefix. 
\end{Deff}

\noindent           The levels of the analytical hierarchy are distinguished by considering the number of alternations in the reduced prefix. 

\begin{Deff}
For $n>0$, a $\Si^1_n$-prefix is one whose reduced prefix begins with $\exists^1$ and has $n-1$ alternations of quantifiers. 
A $\Si^1_0$-prefix is one whose reduced prefix is empty. 
For $n>0$, a $\Pi^1_n$-prefix is one whose reduced prefix begins with $\fa^1$ and has $n-1$ alternations of quantifiers. 
A $\Pi^1_0$-prefix is one whose reduced prefix is empty. 
\nl A predicate form is a $\Si^1_n$ ($\Pi^1_n$)-form if it has a  $\Si^1_n$ ($\Pi^1_n$)-prefix. 
The class of sets in some $\mathbb{N}^l$ which can be expressed in $\Si^1_n$-form (respectively, $\Pi^1_n$-form) is denoted by 
$\Si^1_n$   (respectively, $\Pi^1_n$). 
\nl The class $\Si^1_0 = \Pi^1_0$ is the class of arithmetical sets. 
\end{Deff}

\noindent          We now recall some well known results about the analytical hierarchy. 

\begin{Pro}
Let $R \subseteq \mathbb{N}^l$ for some integer $l$. Then $R$ is an analytical set iff there is some integer $n\geq 0$ such that 
$R \in \Si^1_n$ or $R \in \Pi^1_n$. 
\end{Pro}

\begin{The} For each integer $n\geq 1$, 
\noindent          
\begin{enumerate}
\ite[(a)] $\Si^1_n\cup \Pi^1_n \subsetneq  \Si^1_{n +1}\cap \Pi^1_{n +1}$.
\ite[(b)] A set $R \subseteq \mathbb{N}^l$ is in the class $\Si^1_n$ iff its 
complement is in the class $\Pi^1_n$. 
\ite[(c)] $\Si^1_n - \Pi^1_n \neq \emptyset$ and $\Pi^1_n - \Si^1_n \neq \emptyset$.
\end{enumerate}
\end{The}

\noindent           Transformations of prefixes  are often used, following the rules given by the next theorem. 

\begin{The}
For any predicate form with the given prefix, an equivalent predicate form with the new one can be obtained, following the 
allowed prefix transformations given below :
\noindent          
\begin{enumerate}
\ite[(a)]  $\ldots \exists^0 \exists^0 \ldots \ra \ldots  \exists^0 \ldots, ~~~~~~ \ldots \fa^0  \fa^0 \ldots  \ra \ldots  \fa^0 \ldots ; $
\ite[(b)]  $\ldots \exists^1 \exists^1 \ldots \ra \ldots  \exists^1 \ldots, ~~~~~~ \ldots \fa^1  \fa^1 \ldots  \ra \ldots  \fa^1 \ldots ;$
\ite[(c)]  $\ldots \exists^0 ~~~\ldots  \ra \ldots  \exists^1 \ldots, ~~~~~~ \ldots  \fa^0 ~~~\ldots  \ra \ldots   \fa^1 \ldots ; $
\ite[(d)]  $\ldots \exists^0 \fa^1 \ldots \ra \ldots \fa^1 \exists^0 \ldots, ~~~~ \ldots \fa^0 \exists^1 \ldots  \ra \ldots \exists^1 \fa^0 \ldots ; $
\end{enumerate}
\end{The}

\noindent          We can now define the notion of 1-reduction and of    $\Si^1_n$-complete (respectively,           $\Pi^1_n$-complete) sets. 
Notice that we give the definition for subsets of  $\mathbb{N}$ but one can  easily extend this definition to the case of  
 subsets of $\mathbb{N}^l$ for some integer $l$. 

\begin{Deff}
Given two sets $A,B \subseteq \mathbb{N}$ we say A is 1-reducible to B and write $A \leq_1 B$
if there exists a total computable injective  function f from      $\mathbb{N}$     to   $\mathbb{N}$     such that  $A = f ^{-1}[B]$. 
\end{Deff}

\begin{Deff}
A set $A \subseteq \mathbb{N}$ is said to be $\Si^1_n$-complete   (respectively,   $\Pi^1_n$-complete)  iff $A$ is a  $\Si^1_n$-set 
 (respectively,   $\Pi^1_n$-set) and for each $\Si^1_n$-set  (respectively,   $\Pi^1_n$-set) $B \subseteq \mathbb{N}$ it holds that 
$B \leq_1 A$. 
\end{Deff}

\noindent          For each integer $n\geq 1$ there exists some $\Si^1_n$-complete set $E_n \subseteq \mathbb{N}$. The complement $E_n^-=\mathbb{N}-E_n$ is 
a $\Pi^1_n$-complete set. These sets are precisely defined in \cite{rog} or \cite{cc}.

\subsection{Borel Hierarchy and Analytic Sets}

\noindent          We assume now the reader to be familiar with basic notions of topology which
may be found in \cite{Moschovakis80,LescowThomas,Kechris94,Staiger97,PerrinPin}.

\bigskip \noindent                   There is a natural metric on the set $\Sio$ of  infinite words 
over a finite alphabet 
$\Si$ containing at least two letters which is called the {\it prefix metric} and defined as follows. For $u, v \in \Sio$ and 
$u\neq v$ let $\delta(u, v)=2^{-l_{\mathrm{pref}(u,v)}}$ where $l_{\mathrm{pref}(u,v)}$ 
 is the first integer $n$
such that the $(n+1)^{st}$ letter of $u$ is different from the $(n+1)^{st}$ letter of $v$. 
This metric induces on $\Sio$ the usual  Cantor topology for which {\it open subsets} of 
$\Sio$ are in the form $W.\Si^\om$, where $W\subseteq \Sis$.
A set $L\subseteq \Si^\om$ is a {\it closed set} iff its complement $\Si^\om - L$ 
is an open set.
Now let define  the {\it Borel Hierarchy} of subsets of $\Si^\om$:

\begin{Deff}
For a non-null countable ordinal $\alpha$, the classes ${\bf \Si}^0_\alpha$
 and ${\bf \Pi}^0_\alpha$ of the Borel Hierarchy on the topological space $\Si^\om$ 
are defined as follows:
\nl ${\bf \Si}^0_1$ is the class of open subsets of $\Si^\om$, 
 ${\bf \Pi}^0_1$ is the class of closed subsets of $\Si^\om$, 
\nl and for any countable ordinal $\alpha \geq 2$: 
\nl ${\bf \Si}^0_\alpha$ is the class of countable unions of subsets of $\Si^\om$ in 
$\bigcup_{\gamma <\alpha}{\bf \Pi}^0_\gamma$.
 \nl ${\bf \Pi}^0_\alpha$ is the class of countable intersections of subsets of $\Si^\om$ in 
$\bigcup_{\gamma <\alpha}{\bf \Si}^0_\gamma$.
\end{Deff}

\noindent          For 
a countable ordinal $\alpha$,  a subset of $\Si^\om$ is a Borel set of {\it rank} $\alpha$ iff 
it is in ${\bf \Si}^0_{\alpha}\cup {\bf \Pi}^0_{\alpha}$ but not in 
$\bigcup_{\gamma <\alpha}({\bf \Si}^0_\gamma \cup {\bf \Pi}^0_\gamma)$.

\bigskip \noindent                      
There are also some subsets of $\Si^\om$ which are not Borel. 
Indeed there exists another hierarchy beyond the Borel hierarchy, which is called the 
projective hierarchy and which is obtained from  the Borel hierarchy by 
successive applications of operations of projection and complementation.
The first level of the projective hierarchy is formed by the class of {\it analytic sets} and the class of {\it co-analytic sets} which are complements of 
analytic sets.  
In particular 
the class of Borel subsets of $\Si^\om$ is strictly included into 
the class  ${\bf \Si}^1_1$ of {\it analytic sets} which are 
obtained by projections of Borel sets. 

\begin{Deff} 
A subset $A$ of  $\Si^\om$ is in the class ${\bf \Si}^1_1$ of {\bf analytic} sets
iff there exist a finite set $Y$ and a Borel subset $B$ of $(\Si \times Y)^\om$ 
such that $[ x \in A \lra \exists y \in Y^\om $ ~ $(x, y) \in B  ]$, where $(x, y)$
is the infinite word over the alphabet $\Si \times Y$ such that
$(x, y)(i)=(x(i),y(i))$ for each  integer $i\geq 1$.
\end{Deff} 
 
\noindent           We now define completeness with regard to reduction by continuous functions. 
For a countable ordinal  $\alpha\geq 1$, a set $F\subseteq \Si^\om$ is said to be 
a ${\bf \Si}^0_\alpha$  
(respectively,  ${\bf \Pi}^0_\alpha$, ${\bf \Si}^1_1$)-{\it complete set} 
iff for any set $E\subseteq Y^\om$  (with $Y$ a finite alphabet): 
 $E\in {\bf \Si}^0_\alpha$ (respectively,  $E\in {\bf \Pi}^0_\alpha$,  $E\in {\bf \Si}^1_1$) 
iff there exists a continuous function $f: Y^\om \ra \Si^\om$ such that $E = f^{-1}(F)$. 
 ${\bf \Si}^0_n$
 (respectively ${\bf \Pi}^0_n$)-complete sets, with $n$ an integer $\geq 1$, 
 are thoroughly characterized in \cite{Staiger86a}.  

\bigskip \noindent                   In particular  $\mathcal{R}=(0^\star.1)^\om$  
is a well known example of a 
${\bf \Pi}^0_2 $-complete subset of $\{0, 1\}^\om$. It is the set of 
$\om$-words over $\{0, 1\}$ having infinitely many occurrences of the letter $1$. 
Its  complement 
$\{0, 1\}^\om - (0^\star.1)^\om$ is a 
${\bf \Si}^0_2 $-complete subset of $\{0, 1\}^\om$.

\bigskip \noindent                   We recall now the definition of the  arithmetical hierarchy of  \ol s which form the effective analogue to the 
hierarchy of Borel sets of finite ranks. 
\nl Let $X$ be a finite alphabet. An \ol~ $L\subseteq X^\om$  belongs to the class 
$\Si_n$ if and only if there exists a recursive relation 
$R_L\subseteq (\mathbb{N})^{n-1}\times X^\star$  such that
$$L = \{\sigma \in X^\om \mid \exists a_1\ldots Q_na_n  \quad (a_1,\ldots , a_{n-1}, 
\sigma[a_n+1])\in R_L \}$$

\noindent          where $Q_i$ is one of the quantifiers $\fa$ or $\exists$ 
(not necessarily in an alternating order). An \ol~ $L\subseteq X^\om$  belongs to the class 
$\Pi_n$ if and only if its complement $X^\om - L$  belongs to the class 
$\Si_n$.  The inclusion relations that hold  between the classes $\Si_n$ and $\Pi_n$ are 
the same as for the corresponding classes of the Borel hierarchy. 
 The classes $\Si_n$ and $\Pi_n$ are  included in the respective classes 
${\bf \Si_n^0}$ and ${\bf \Si_n^0}$ of the Borel hierarchy, and cardinality arguments suffice to show that these inclusions are strict. 

\bigskip \noindent                    As in the case of the Borel hierarchy, projections of arithmetical sets 
(of the second $\Pi$-class) lead 
beyond the arithmetical hierarchy, to the analytical hierarchy of \ol s. The first class 
of this hierarchy is the (lightface) class $\Si^1_1$ of {\it effective analytic sets} 
 which are obtained by projection of arithmetical sets.
An \ol~ $L\subseteq X^\om$  belongs to the class 
$\Si_1^1$ if and only if there exists a recursive relation 
$R_L\subseteq \mathbb{N}\times \{0, 1\}^\star \times X^\star$  such that:

$$L = \{\sigma \in X^\om  \mid \exists \tau (\tau\in \{0, 1\}^\om \wedge \fa n \exists m 
 ( (n, \tau[m], \sigma[m]) \in R_L )) \}$$

\noindent          Then an \ol~ $L\subseteq X^\om$  is in the class $\Si_1^1$ iff it is the projection 
of an \ol~ over the alphabet $X\times \{0, 1\}$ which is in the class $\Pi_2$.  The (lightface)  class $\Pi_1^1$ of  {\it effective co-analytic sets} 
 is simply the class of complements of effective analytic sets. We denote as usual $\Delta_1^1 = \Si^1_1 \cap \Pi_1^1$. 
\nl Recall that an \ol~ $L\subseteq X^\om$ is in the class $\Si_1^1$
iff it is accepted by a non deterministic Turing machine (reading $\om$-words)
with a   B\"uchi or Muller acceptance condition  \cite{CG78b,Staiger97}.

\bigskip \noindent                    For $\Ga$ a finite alphabet having at least two letters, the 
set $\Ga^{\om \times \om}$ of functions  from $\om \times \om$ into $\Ga$ 
is usually equipped with the product topology  of the discrete 
topology on $\Ga$. 
This topology may be defined 
by the following distance $d$. Let $x$ and $y$  in $\Ga^{\om  \times \om}$ 
such that $x\neq y$, then  
$$ d(x, y)=\frac{1}{2^n} ~~~~~~~\mbox{   where  }$$
$$n=min\{p\geq 0 \mid  \exists (i, j) ~~ x(i, j)\neq y(i, j) \mbox{ and } i+j=p\}.$$

\noindent           Then the topological space $\Ga^{\om \times \om}$ is homeomorphic to the 
topological space $\Ga^{\om}$, equipped with the Cantor topology.  
Borel subsets  of   $\Ga^{\om  \times \om}$ are defined from open 
subsets  as in the case of the topological space $\Ga^\om$. 
Analytic  subsets  of   $\Ga^{\om  \times \om}$ are obtained as projections on 
$\Ga^{\om \times \om}$
 of Borel subsets of the product space  $\Ga^{\om \times \om} \times \Ga^\om$.
\nl  The set $\Si^{\om, \om}$ of $\om$-pictures over $\Si$, 
considered  a topological subspace of $\hat{\Si}^{\om \times \om}$, 
 is easily seen to be homeomorphic to the topological space $\Si^{\om \times \om}$, 
via the mapping 
$\varphi: \Si^{\om, \om} \ra \Si^{\om \times \om}$ 
defined by $\varphi(p)(i, j)=p(i+1, j+1)$ for all 
$p\in \Si^{\om, \om}$ and $i, j \in \om$.

 \subsection{Some Results  of Set Theory}

\noindent          We now recall some basic notions of set theory 
which will be useful in the sequel, and which are exposed in any  textbook on set theory, like \cite{Jech}.

\bigskip \noindent                    The usual axiomatic system {\bf ZFC} is 
Zermelo-Fraenkel system {\bf ZF}   plus the axiom of choice {\bf AC}. 
 A model ({\bf V}, $\in)$ of  the axiomatic system {\bf ZFC} is a collection  {\bf V} of sets,  equipped with 
the membership relation $\in$, where ``$x \in y$" means that the set $x$ is an element of the set $y$, which satisfies the axioms of  {\bf ZFC}.  
We shall often say `` the model {\bf V}"
instead of  ``the model  ({\bf V}, $\in)$". 

\bigskip \noindent                   The axioms of {\bf ZFC} express some  natural facts that we consider to hold in the universe of sets. For instance a natural fact is that 
two sets $x$ and $y$ are equal iff they have the same elements. 
This is expressed by the {\it Axiom of Extensionality}. 
  Another natural axiom is the {\it Pairing Axiom}   which states that for all sets $x$ and $y$ there exists a  set $z=\{x, y\}$ 
whose elements are $x$ and $y$. 
Similarly the {\it Powerset Axiom} states the existence of the set of subsets of a set $x$. 
We refer the reader to any textbook on set theory, like \cite{Jech},  for an exposition of the other axioms of {\bf ZFC}.

\bigskip \noindent                    The infinite cardinals are usually denoted by
$\aleph_0, \aleph_1, \aleph_2, \ldots , \aleph_\alpha, \ldots$
The cardinal $\aleph_\alpha$ is also denoted by $\om_\alpha$,
as usual when it is considered an ordinal.

\bigskip \noindent                   The continuum hypothesis {\bf CH}  says that the first uncountable cardinal $\aleph_1$ is equal to $2^{\aleph_0}$ which is the cardinal of the 
continuum. G\"odel and Cohen  proved that the continuum hypothesis {\bf CH} is independent from the axiomatic system {\bf ZFC}: providing 
{\bf ZFC} is consistent, there exist some models of {\bf ZFC + CH} 
and also some models of {\bf ZFC + $\neg$ CH}, where {\bf $\neg$ CH} denotes the negation of the 
continuum hypothesis, \cite{Jech}. 

\bigskip \noindent                   
Let ${\bf ON}$ be the class of all ordinals. Recall that an ordinal $\alpha$ is said to be a successor ordinal iff there exists an ordinal $\beta$ such that 
$\alpha=\beta + 1$; otherwise the ordinal $\alpha$ is said to be a limit ordinal and in that case 
$\alpha ={\rm sup} \{ \beta \in {\bf ON}\mid \beta < \alpha \}$.

\bigskip \noindent                    The  class ${\bf L}$ of  {\it constructible sets} in a model {\bf V} of {\bf ZF} is defined by 
$${\bf L} = \bigcup_{\alpha \in {\bf ON}} {\bf L}(\alpha) $$
\noindent          where the sets ${\bf L}(\alpha) $ are constructed  by induction as follows: 

\begin{enumerate}
\ite ${\bf L}(0) =\emptyset$
\ite ${\bf L}(\alpha) = \bigcup_{\beta <  \alpha} {\bf L}(\beta) $, for $\alpha$ a limit ordinal, and 
\ite ${\bf L}(\alpha + 1) $ is the set of subsets of ${\bf L}(\alpha) $ which are definable from a finite number of elements of ${\bf L}(\alpha) $
by a first-order formula relativized to ${\bf L}(\alpha) $. 
\end{enumerate}

\noindent           If  {\bf V} is  a model of {\bf ZF} and ${\bf L}$ is  the class of  {\it constructible sets} of   {\bf V}, then the class  ${\bf L}$     forms a model of  
{\bf ZFC + CH}.
Notice that the axiom ({\bf V=L}) means ``every set is constructible"  and that it is consistent with {\bf ZFC}.  

\bigskip \noindent                   Consider now a model {\bf V} of  the axiomatic system {\bf ZFC} and the class of constructible sets ${\bf L} \subseteq {\bf V}$ which forms another 
model of  {\bf ZFC}.  It is known that 
the ordinals of {\bf L} are also the ordinals of  {\bf V}. But the cardinals  in  {\bf V}  may be different from the cardinals in {\bf L}. 

\bigskip \noindent                    In particular,  the first uncountable cardinal in {\bf L} is denoted 
 $\aleph_1^{\bf L}$. It is in fact an ordinal of {\bf V} which is denoted $\om_1^{\bf L}$. 
  It is known that this ordinal satisfies the inequality 
$\om_1^{\bf L} \leq \om_1$.  In a model {\bf V} of  the axiomatic system {\bf ZFC + V=L} the equality $\om_1^{\bf L} = \om_1$ holds. But in 
some other models of {\bf ZFC} the inequality may be strict and then $\om_1^{\bf L} < \om_1$. This is explained in \cite[page 202]{Jech}: one can start 
 from a model 
{\bf V} of {\bf ZFC + V=L} and construct by  forcing  a generic extension {\bf V[G]} in which the cardinals $\om$ and $\om_1$ are 
collapsed; in this extension the inequality $\om_1^{\bf L} < \om_1$ holds.

\bigskip \noindent                   We now recall the notion of  a perfect set. 

\begin{Deff} 
Let $P \subseteq \Sio$, where $\Si$ is a finite alphabet having at least two letters. The set 
 $P$ is said to be a perfect subset of $\Sio$ if and only if :  
\nl  (1) $P$ is a non-empty closed set,  and 
\nl (2) for every $x\in P$ and every open set $U$ containing $x$ there is an element 
$y \in P\cap U$ such that $x\neq y$. 
\end{Deff}

\noindent          So a perfect subset of $\Sio$ is a non-empty closed set which has no isolated points. It is well known that a  perfect subset of 
$\Sio$   has cardinality $2^{\aleph_0}$, see \cite[page 66]{Moschovakis80}. 

\bigskip \noindent                   We now  recall   the notion of  thin  subset of $\Sio$. 

 \begin{Deff}
A set $X \subseteq \Sio$ is said to be thin iff it  contains no perfect subset. 
\end{Deff}

\noindent          The following  
important  result was proved by Kechris \cite{Kechris75} and independently by Guaspari \cite{Guaspari}  and  Sacks  \cite{Sacks}. 

\begin{The}[see \cite{Moschovakis80} page 247] ({\bf ZFC})  Let $\Si$ be a finite alphabet having at least two letters. 
There exists a thin $\Pi_1^1$-set $\mathcal{C}_1( \Sio) \subseteq  \Sio$ which contains every thin,  $\Pi_1^1$-subset of $\Sio$. 
It is called the  largest thin $\Pi_1^1$-set  in $\Sio$.    
\end{The}

\noindent          An important fact is that the cardinality of the largest thin $\Pi_1^1$-set in $\Sio$  depends on the model of {\bf ZFC}. 
The following  result on the cardinality of $\mathcal{C}_1( \Sio) $, was proved by Kechris and  independently by Guaspari and Sacks, see also 
\cite[page 171]{Kanamori}. 

\begin{The}
({\bf ZFC})   The cardinal  of the  largest thin $\Pi_1^1$-set in  $\Sio$ is equal to the cardinal of  $\om_1^{\bf L}$. 
\end{The}

\noindent          This means that in a given model {\bf V} of {\bf ZFC} the cardinal  of the  largest thin $\Pi_1^1$-set in  $\Sio$ is equal to the cardinal 
{\it in {\bf V}} of the ordinal $\om_1^{\bf L}$ which plays the role of the cardinal $\aleph_1$ in the inner model {\bf L}  of constructible sets of {\bf V}.

\bigskip \noindent                   We can now state the following result which will be useful in the sequel. 

\begin{Cor}\label{cor1}
\noindent          
\begin{itemize}
\item[(a)]   There is a model ${\bf V}_1$ of {\bf ZFC} in which the largest thin $\Pi_1^1$-set in  $\Sio$ has cardinal $\aleph_1$,  where  
 $\aleph_1=2^{\aleph_0}$. 
\item[(b)]     There is a model ${\bf V}_2$ of {\bf ZFC} in which the largest thin $\Pi_1^1$-set in  $\Sio$ has cardinal $\aleph_0$, 
i.e.  is countably infinite. 
\item[(c)]  There is a model ${\bf V}_3$ of {\bf ZFC} in which the largest thin $\Pi_1^1$-set in  $\Sio$ has cardinal $\aleph_1$, where  
$\aleph_0 < \aleph_1 < 2^{\aleph_0}$. 
\end{itemize}
\end{Cor}

\proof  (a).  In the model {\bf L},  the cardinal of the  largest thin $\Pi_1^1$-set in  $\Sio$ is equal to the cardinal of  $\om_1$. 
Moreover the continuum hypothesis is satisfied thus $2^{\aleph_0}=\aleph_1$. 
\nl  Thus the largest thin $\Pi_1^1$-set in  $\Sio$ has the cardinality $2^{\aleph_0}=\aleph_1$.

\bigskip \noindent                   (b).  Let {\bf V} be a model of  ({\bf ZFC} + $\om_1^{\bf L} < \om_1$).  In this model $\om_1$ is the first uncountable ordinal. Thus 
$\om_1^{\bf L} < \om_1$ implies that $\om_1^{\bf L}$ is a countable ordinal in {\bf V}. Its cardinal is $\aleph_0$ and it is also the cardinal of 
the  largest thin $\Pi_1^1$-set in  $\Sio$.  

\bigskip \noindent                   (c).  It suffices to show that there is a model ${\bf V}_3$ of {\bf ZFC} in which $\om_1^{\bf L} = \om_1$ and $\aleph_1 < 2^{\aleph_0}$. 
Such a model can be constructed by Cohen's forcing. We can start  from a model 
{\bf V} of {\bf ZFC + V=L} (in which   $\om_1^{\bf L} = \om_1$)        and construct by  forcing  a generic extension {\bf V[G]} 
in which $\aleph_2$ subsets of $\om$ are added. 
 Notice that the cardinals are preserved under this extension (see 
\cite[page 219]{Jech}) and that the constructible sets of  {\bf V[G]} are also  the   constructible sets of    {\bf V}. 
 Thus in the new model {\bf V[G]} we still have 
$\om_1^{\bf L} = \om_1$ but now $\aleph_1 < 2^{\aleph_0}$. 
\ep

\section{Decision  Problems}

\noindent           We  now  study decision problems for recognizable languages of infinite pictures. 
We gave  in  \cite{Fink-tilings}  the exact degree of several natural decision problems. We   first recall some of these results. 

\bigskip \noindent                     Castro and Cucker proved in \cite{cc} that the non-emptiness problem and the infiniteness problem for $\om$-languages 
of Turing machines are both
$\Si_1^1$-complete. We  easily inferred  from this result  a similar result for recognizable languages of infinite pictures. 

\bigskip \noindent                   From now on we shall  denote by ${\cal T}_z$ the non deterministic tiling system of index $z$, 
(accepting pictures over $\Si=\{a, b\}$),
 equipped with a B\"uchi acceptance condition. 

\begin{The}[\cite{Fink-tilings}]\label{E-I}
The  non-emptiness  problem and the infiniteness problem for B\"uchi-recognizable languages of infinite pictures 
 are $\Si_1^1$-complete, i.e. : 
\begin{enumerate}
\ite  $\{ z \in \mathbb{N} \mid  L^B({\cal T}_z) \neq  \emptyset  \}$ is  $\Si_1^1$-complete.  
\ite  $\{ z \in \mathbb{N} \mid  L^B({\cal T}_z)  \mbox{ is infinite }   \}$ is  $\Si_1^1$-complete.
\end{enumerate}
\end{The}
  
\noindent          In a similar way, the universality  problem and the inclusion and the equivalence problems,  for $\om$-languages 
of Turing machines, have been proved to be $\Pi_2^1$-complete by Castro and Cucker in  \cite{cc}, and we used these results to prove 
 the following results in \cite{Fink-tilings}. 

\begin{The}[\cite{Fink-tilings}]\label{U}
The universality  problem for B\"uchi-recognizable languages of infinite pictures 
 is $\Pi_2^1$-complete, i.e. : ~~
 $$\{ z \in \mathbb{N} \mid  L^B({\cal T}_z) =  \Si^{\om, \om} \}\mbox{ is  } \Pi_2^1 \mbox{-complete.}$$ 
\end{The}

\begin{The}[\cite{Fink-tilings}]\label{inclusion-equivalence}
\noindent           The inclusion and the equivalence problems for B\"uchi-recognizable languages of infinite pictures are 
  $\Pi_2^1$-complete, i.e. : 
\begin{enumerate}
\ite $\{ (y, z) \in \mathbb{N}^2  \mid  L^B({\cal T}_y) \subseteq L^B({\cal T}_z)  \}$ is  $\Pi_2^1$-complete. 
\ite $\{ (y, z) \in \mathbb{N}^2  \mid  L^B({\cal T}_y) = L^B({\cal T}_z)  \}$ is  $\Pi_2^1$-complete. 
\end{enumerate}
\end{The}

\bigskip \noindent                    The class of      B\"uchi-recognizable  languages of   infinite pictures 
is not closed under complement \cite{ATW02}. 
Thus the following question naturally arises: ``can we decide whether the complement of a 
 B\"uchi-recognizable   language of  infinite pictures is B\"uchi-recognizable?". And what is the exact complexity of this decision problem, called the 
complementability  problem. 

\bigskip \noindent                   Another classical problem is the determinizability problem: ``can we decide whether a given  recognizable   language of  infinite pictures is recognized by a 
deterministic tiling system?". 
\nl  Recall  that
 a tiling system is called deterministic if on any picture 
it allows at most one tile covering the origin, the state assigned to position 
$(i+1, j+1)$ is uniquely determined by the states at positions $(i, j), (i+1, j), (i, j+1)$ 
and the states at the border positions $(0, j+1)$ and $(i+1, 0)$ are determined by the state 
$(0, j)$, respectively $(i, 0)$, \cite{ATW02}. 
\nl As remarked in \cite{ATW02}, the hierarchy proofs of the classical 
Landweber hierarchy defined using deterministic $\om$-automata ``carry over without essential 
changes to pictures". In particular, a language of $\om$-pictures which 
is B\"uchi-recognized by a \de  tiling system is a ${\bf \Pi}^0_2$-set and a language of $\om$-pictures which 
is Muller-recognized by a \de  tiling system is a boolean combination of ${\bf \Pi}^0_2$-sets, hence a ${\bf \Delta}^0_3$-set. 

\bigskip \noindent                   These topological properties have been used in \cite{Fink-tilings}, along with a dichotomy property, to prove the following results. 

\begin{The}[\cite{Fink-tilings}]
\noindent           The determinizability problem and the complementability  problem for B\"uchi-recognizable languages of infinite pictures 
are   $\Pi_2^1$-complete,  i.e. : 
\begin{enumerate}
\ite 
$\{  z \in \mathbb{N}  \mid  
L^B({\cal T}_z) \mbox{ is B\"uchi-recognizable  by a deterministic tiling system}\}$ is $\Pi_2^1$-complete. 
\ite
$\{  z \in \mathbb{N}  \mid 
 L^B({\cal T}_z) \mbox{ is Muller-recognizable  by a deterministic tiling system }\}$ is $\Pi_2^1$-complete. 
\ite
 $\{  z \in \mathbb{N}  \mid  \exists y ~~ \Si^{\om, \om} - L^B({\cal T}_z) =  L^B({\cal T}_y) \}$ is $\Pi_2^1$-complete. 
 \end{enumerate}
\end{The}

\noindent         
We already mentioned that we used some results of Castro and Cucker in the proof of the above results. 
Castro and Cucker  studied degrees of decision problems for $\om$-languages accepted by Turing machines and 
proved that many of them are highly undecidable, \cite{cc}. We are going  to use again some of their results to prove here new results about 
B\"uchi-recognizable languages of infinite pictures. 

\bigskip \noindent                     We firstly recall  the notion of acceptance of infinite words by Turing machines considered by Castro and Cucker in \cite{cc}. 

\begin{Deff}
A non deterministic Turing machine ${\cal M}$ is a $5$-tuple ${\cal M}=(Q, \Si, \Ga, \delta, q_0)$, where $Q$ is a finite set of states, 
$\Si$ is a finite input alphabet, $\Ga$ is a finite tape alphabet satisfying $\Si  \subseteq \Ga$, $q_0$ is the initial state, 
and $\delta$ is a mapping from $Q \times \Ga$ to subsets of $Q \times \Ga \times \{L, R, S\}$. A configuration of ${\cal M}$ is a $3$-tuple 
$(q, \sigma, i)$, where $q\in Q$, $\sigma \in \Ga^\om$ and $i\in \mathbb{N}$. An infinite sequence of configurations $r=(q_i, \alpha_i, j_i)_{i\geq 1}$
is called a run of ${\cal M}$ on $w\in \Sio$ iff: 
\begin{enumerate}
\ite[(a)] $(q_1, \alpha_1, j_1)=(q_0, w, 1)$, and 
\ite[(b)] for each $i\geq 1$, $(q_i, \alpha_i, j_i) \vdash (q_{i+1}, \alpha_{i+1}, j_{i+1})$, 
\end{enumerate}
\noindent          where $\vdash$ is the transition relation of ${\cal M}$ defined as usual. The run $r$ is said to be complete if 
 the limsup of the head positions is infinity, i.e. if $(\fa n \geq 1) (\exists k \geq 1) (j_k \geq n)$. 
The run $r$ is said to be oscillating if the liminf of the head positions is bounded, i.e. if $(\exists k \geq 1) (\fa n \geq 1) (\exists m \geq n) ( j_m=k)$. 

\end{Deff}

\begin{Deff}
Let ${\cal M}=(Q, \Si, \Ga, \delta, q_0)$ be a non deterministic Turing machine   and $F \subseteq Q$. The $\om$-language accepted by $({\cal M}, F)$ is 
the set of $\om$-words $ \sigma \in \Sio$ such that there exists a complete non oscillating run $ r=(q_i, \alpha_i, j_i)_{i\geq 1}$
 of ${\cal M}$  on  $\sigma$ such that, for all $ i, q_i \in F.$
\end{Deff}

\noindent          The above acceptance condition is denoted $1'$-acceptance in  \cite{CG78b}. Another usual acceptance condition is the now called B\"uchi 
acceptance condition which is also denoted $2$-acceptance in  \cite{CG78b}. 
We now  recall its definition. 

\begin{Deff}
Let ${\cal M}=(Q, \Si, \Ga, \delta, q_0)$ be a non deterministic Turing machine   
and $F \subseteq Q$. The $\om$-language B\"uchi accepted by $({\cal M}, F)$ is 
the set of $\om$-words $ \sigma \in \Sio$ such that there exists a complete non oscillating run $ r=(q_i, \alpha_i, j_i)_{i\geq 1}$
 of ${\cal M}$  on  $\sigma$ and  infinitely many integers $i$ such that $q_i \in F.$
\end{Deff}

\noindent          Recall that Cohen and Gold proved in \cite[Theorem 8.6]{CG78b} that one can effectively construct, from a given non deterministic Turing machine, 
another equivalent non deterministic Turing machine, equipped with the same kind of  acceptance condition, and in which every run is complete non oscillating. 
 Cohen and Gold proved also in \cite[Theorem 8.2]{CG78b} that an $\om$-language is accepted by a non deterministic Turing machine with 
$1'$-acceptance condition iff it is accepted by a non deterministic Turing machine with B\"uchi acceptance condition.

\bigskip \noindent                   From now on, we shall  denote ${\cal M}_z$ the non deterministic Turing machine of index $z$, (accepting words over $\Si=\{a, b\}$), 
equipped with a $1'$-acceptance condition. 

\bigskip \noindent                   An important notion in automata theory is the notion of ambiguity. It can be defined also in the context of acceptance by tiling systems, see 
\cite{AGMR} for the case of finite pictures. 

\begin{Deff}
Let ${\cal A}$=$(Q, \Si, \Delta)$ 
be a tiling system, and $F\subseteq Q$.  The B\"uchi tiling system $({\cal A}, F)$ is unambiguous iff every $\om$-picture $p \in  \Si^{\om, \om}$ has at most 
an accepting run by $({\cal A}, F)$. 
\end{Deff}

\begin{Deff}
A B\"uchi recognizable  language $L \subseteq \Si^{\om, \om}$ is unambiguous iff 
there exists an unambiguous B\"uchi tiling system $({\cal A}, F)$ 
such that $L = L( {\cal A}, F )$. Otherwise the language $L $ is said to be inherently ambiguous. 
\end{Deff}

\noindent          We can now prove the following result, which is very similar to a corresponding result for recognizable tree languages proved in \cite{FS-tree}. 

\begin{Pro}
Let $L \subseteq \Si^{\om, \om}$ be an unambiguous B\"uchi recognizable  language of  infinite pictures. Then $L$ is a Borel subset of $\Si^{\om, \om}$. 
\end{Pro}

\proof Let $L \subseteq  \Si^{\om, \om}$ be a   language  accepted by an unambiguous  
B\"uchi tiling system $({\cal A}, F)$, where  ${\cal A}$=$(Q, \Si, \Delta)$, and let $R \subseteq (\hat{\Si}\times Q)^{\om\times \om}$ be defined  by: 
$$R= \{ (p, \rho) \mid p\in  \Si^{\om, \om} \mbox{ and }   \rho \in   \mbox{ is an accepting run of }({\cal A}, F) \mbox{ on  the picture } p \}.$$
\noindent          The set $R$ is easily seen to be  a  ${\bf \Pi}_2^0$-subset of $(\hat{\Si}\times Q)^{\om\times \om}$. 

\bigskip \noindent                    Consider now the projection $\mathrm{PROJ}_{\hat{\Si}^{\om\times \om}} : ~ \hat{\Si}^{\om\times \om} \times Q^{\om\times \om}
 \ra \hat{\Si}^{\om\times \om}$ defined by 
$\mathrm{PROJ}_{\hat{\Si}^{\om\times \om}}((p, \rho)) = p$ for all $(p, \rho) \in  \hat{\Si}^{\om\times \om} \times Q^{\om\times \om}$. This projection is a continuous function 
and it is {\it injective} on the Borel set $R$ because the B\"uchi tiling system $({\cal A}, F)$ is unambiguous. 
Hence, by a Theorem of Lusin and Souslin, see \cite[Theorem 15.1 page 89]{Kechris94}, the injective 
image of $R$ by the continuous function $\mathrm{PROJ}_{\hat{\Si}^{\om\times \om}}$ is  Borel. 
Thus the  language $L=\mathrm{PROJ}_{\hat{\Si}^{\om\times \om}}(R)$ is a Borel subset of $\hat{\Si}^{\om\times \om}$.  But 
$\Si^{\om, \om}$ is a closed subset of $\hat{\Si}^{\om\times \om}$ and $L \subseteq \Si^{\om, \om}$. 
Thus  $L$ is also a Borel subset of $\Si^{\om, \om}$. 
\ep

\begin{Cor}
There exist some inherently ambiguous B\"uchi-recognizable languages  of infinite pictures. 
\end{Cor}

\proof The result follows directly from the above proposition because we know that there exist some
B\"uchi-recognizable languages  of infinite pictures which are not Borel sets, see \cite{Fin04,Fink-tilings}. 
\ep 

\bigskip \noindent                   We can now state that the unambiguity problem for recognizable  language of  infinite pictures is $\Pi_2^1$-complete.

\begin{The}\label{unamb}
\noindent           The unambiguity  problem      for  recognizable  languages of  infinite pictures  is 
 $\Pi_2^1$-complete, i.e. :  
 $$ \{  z \in \mathbb{N}  \mid  L^B({\cal T}_z) \mbox{ is non ambiguous  }\}\mbox{  is } \Pi_2^1\mbox{-complete.} $$
\end{The}

\proof   To prove that the  unambiguity  problem      for  recognizable  language of  infinite pictures  is in the class 
 $\Pi_2^1$, we reason as in the case of the  unambiguity  problem      for  $\om$-languages accepted by $1$-counter or $2$-tape 
automata, see  \cite{Fin-HI}. 

\bigskip \noindent                   Notice first, as in \cite{Fink-tilings},  that, using a recursive bijection $b : (\mathbb{N}-\{0\})^2  \ra \mathbb{N}-\{0\}$,  
one can associate with  each  $\om$-word $\sigma \in \Sio$ a unique  $\om$-picture $p^\sigma \in \Si^{\om, \om}$ which is simply defined 
by $p^\sigma(i, j)=\sigma(b(i, j))$ for all integers $i, j \geq 1$. And we 
can identify a run $\rho \in Q^{\om \times \om}$ with an element of $Q^\om$ 
and finally with a coding of this element over the alphabet $\{0, 1\}$. So the run  $\rho$ can be identified 
with its code $\bar{\rho} \in  \{0, 1\}^\om$. 

\bigskip \noindent                    If a tiling system  ${\cal A}$=$(Q, \Si, \Delta)$ is equipped with a set of  accepting states $F \subseteq Q$, then  
for $\sigma \in \Sio$ and $\rho \in  \{0, 1\}^\om$, ~~  ``$\rho$ is a B\"uchi accepting run of $({\cal A}, F)$ over the $\om$-picture $p^\sigma$"
can be expressed by an arithmetical formula, see   \cite[Section 2.4]{ATW02} and \cite{Fink-tilings}. 

\bigskip \noindent                   We can now first express 
``${\cal T}_z$ is non ambiguous"  by : 
$$``\fa \sigma \in  \Si^\om  \fa \rho, \rho' \in \{0, 1\}^\om  
[(  \rho \mbox{ and } \rho' \mbox{ are accepting runs of  } {\cal T}_z  \mbox{ on } p^{\sigma}) \ra  \rho=\rho']"$$
\noindent          which is a $\Pi_1^1$-formula. 
 Then `` $L^B({\cal T}_z) $ is non ambiguous" can be expressed by the following formula:  
 ``$\exists y [ L^B({\cal T}_z) = L^B({\cal T}_y)  \mbox{ and }  {\cal T}_y   \mbox{ is non ambiguous}]$". This is a $\Pi_2^1$-formula because 
$L^B({\cal T}_z) = L^B({\cal T}_y) $ can be expressed by a $\Pi_2^1$-formula, and the quantification  $\exists y$ is of type $0$. 
Thus the set $\{  z \in \mathbb{N}  \mid  L^B({\cal T}_z) \mbox{ is non ambiguous  }\}$  is  a $\Pi_2^1$-set.

\bigskip \noindent                   To prove the completeness part of the theorem, 
we shall use the following dichotomy result proved  in \cite[proof of Theorem 5.11]{Fink-tilings}. 
\nl There exists an 
  injective computable function $H \circ \theta$  from $\mathbb{N}$ into $\mathbb{N}$ such that:
\nl {\bf First case:}  If $L({\cal M}_z)= \Si^\om$ then $L^B({\cal T}_{H \circ \theta (z) })= \Si^{\om, \om}$. 
\nl {\bf Second case:} If $L({\cal M}_z) \neq  \Si^\om$ then $L^B({\cal T}_{H \circ \theta (z) })$ is not a Borel set. 

\bigskip \noindent                   In the first case $L^B({\cal T}_{H \circ \theta (z) })= \Si^{\om, \om}$ is obviously an unambiguous language. And in the second case 
the language $L^B({\cal T}_{H \circ \theta (z) })$ cannot be unambiguous because it is not a Borel subset of $\Si^{\om, \om}$.  
Thus, using the reduction $H \circ \theta $, we see that : 
 $$\{ z \in \mathbb{N} \mid  L({\cal M}_z) = \Si^\om \} \leq_1 \{  z \in \mathbb{N}  \mid  L^B({\cal T}_z) \mbox{ is non ambiguous  }\}$$
\noindent          and the result follows from the  $\Pi_2^1$-completeness of the universality problem for $\om$-languages of Turing machines proved by Castro and Cucker in 
\cite{cc}. 
\ep 

\bigskip \noindent                   Notice that the same dichotomy result above 
 with the reduction $H \circ \theta $ was used in \cite{Fink-tilings} to prove that topological properties of 
recognizable languages of infinite pictures are actually highly undecidable. 

\begin{The}[\cite{Fink-tilings}]
\noindent          Let $\alpha$ be a non-null countable ordinal. Then  
\begin{enumerate}
\ite $ \{  z \in \mathbb{N}  \mid   L^B({\cal T}_z) \mbox{ is in the Borel class } {\bf \Si}^0_\alpha \}$ is  $\Pi_2^1$-hard. 
\ite  $ \{  z \in \mathbb{N}  \mid   L^B({\cal T}_z)  \mbox{ is in the Borel class } {\bf \Pi}^0_\alpha \}$ is  $\Pi_2^1$-hard. 
\ite  $ \{  z \in \mathbb{N}  \mid  L^B({\cal T}_z)  \mbox{ is a  Borel set } \}$ is  $\Pi_2^1$-hard. 
\end{enumerate}
\end{The}

\noindent          A natural question is to study similar problems by replacing Borel classes by the effective classes of the arithmetical hierarchy. This was not 
studied in \cite{Fink-tilings}, but a similar problem was solved in  \cite{Fin-HI} for  $\om$-languages accepted by $1$-counter or $2$-tape 
B\"uchi automata. We can reason in a similar way for the case of recognizable languages of infinite pictures, and state  the following result. 

\begin{The}
\noindent          Let $n \geq 1$ be an integer. Then 
\begin{enumerate}
\ite $\{  z \in \mathbb{N}  \mid  L^B({\cal T}_z) \mbox{ is in the arithmetical class }  \Si_n \}$ is  $\Pi_2^1$-complete. 
\ite $\{  z \in \mathbb{N}  \mid  L^B({\cal T}_z) \mbox{ is in the arithmetical class }  \Pi_n \}$ is  $\Pi_2^1$-complete. 
\ite $\{  z \in \mathbb{N}  \mid  L^B({\cal T}_z) \mbox{ is a }  \Delta^1_1 \mbox{-set } \}$ is  $\Pi_2^1$-complete. 

\end{enumerate}
\end{The}

\noindent          We do not give the complete proof here. It is actually  very similar to the case of $\om$-languages accepted by $1$-counter or $2$-tape 
B\"uchi automata in \cite{Fin-HI}. 
A key argument, to prove that 
$\{  z \in \mathbb{N}  \mid  L^B({\cal T}_z) \mbox{ is in the arithmetical class }  \Si_n \}$
(respectively, $\{  z \in \mathbb{N}  \mid  L^B({\cal T}_z) \mbox{ is in the arithmetical class }  \Pi_n \}$)  is a $\Pi_2^1$-set,  
 is the   existence of a universal set 
$\mathcal{U}_{ \Si_n}\subseteq \mathbb{N} \times  \Si^{\om, \om}$ 
(respectively,  $\mathcal{U}_{ \Pi_n}\subseteq \mathbb{N} \times  \Si^{\om, \om}$) 
   for the class of 
$\Si_n$-subsets of $\Si^{\om, \om}$,   (respectively,     $\Pi_n$-subsets of $\Si^{\om, \om}$),    \cite[p. 172]{Moschovakis80}. 
Notice also that the completeness part follows easily from  the dichotomy result obtained with the reduction $H \circ \theta $.

\bigskip \noindent                   We now come to cardinality problems. We already know that it is $\Si_1^1$-complete to determine whether a given recognizable language
 of infinite pictures is empty (respectively, infinite).  Recall that every  recognizable language
 of infinite pictures is an analytic set. On the other hand, every analytic set is either countable or has  the cardinality $2^{\aleph_0}$ of the continuum. 
Then some  questions naturally arise. What are the complexities of the following decision problems: ``Is a given   recognizable language 
of infinite pictures countable? Is it countably infinite? Is it uncountable?".  Notice that similar  questions were asked by Castro and Cucker 
in the case of  $\om$-languages of Turing machines and have been solved very recently by Finkel and Lecomte in \cite{FL-TM}. 
We can now state the following result for recognizable languages of infinite pictures. Below $D_2(\Sigma_1^1)$ denotes the class of $2$-differences of 
$\Sigma_1^1$-sets, i.e. the class of sets which are intersections of a $\Sigma_1^1$-set and of a $\Pi_1^1$-set.

\begin{The}
\noindent         
\begin{enumerate}
\ite $\{ z \in \mathbb{N} \mid  L^B({\cal T}_z) \mbox{ is  countable} \}$ is $\Pi_1^1$-complete. 
\ite $\{ z \in \mathbb{N} \mid  L^B({\cal T}_z) \mbox{ is uncountable} \}$ is $\Sigma_1^1$-complete. 
\ite $\{ z \in \mathbb{N} \mid  L^B({\cal T}_z) \mbox{  is countably infinite} \}$ is $D_2(\Sigma_1^1)$-complete. 
\end{enumerate}
\end{The}

\proof (1). We can first prove that $\{ z \in \mathbb{N} \mid  L^B({\cal T}_z) \mbox{ is  countable} \}$ is in the class $\Pi_1^1$ in the same way as 
in the case of $\om$-languages of Turing machines in \cite{FL-TM}. 
\nl We know that a  recognizable  language of  infinite pictures   $L^B({\cal T}_z) $      is a $\Sigma_1^1$-subset of  $\Si^{\om, \om}$. 
But it is known that a $\Sigma_1^1$-subset $L$ of  $\Si^{\om, \om}$    is  countable 
 if and only if for every $x\in L$ the singleton  $\{x\}$ is a $\Delta_1^1$-subset 
of $\Si^{\om, \om}$, see \cite[page 243]{Moschovakis80}.
Then, using a nice coding of $\Delta_1^1$-subsets of $\Si^{\om, \om}$ given in \cite[Theorem 3.3.1]{HKL}, we can  
prove that $\{ z \in \mathbb{N} \mid  L^B({\cal T}_z) \mbox{ is  countable} \}$ is in the class $\Pi_1^1$, see \cite{FL-TM} for more details. 

\bigskip \noindent                   To prove the completeness part of Item (1), we shall use the following two  lemmas proved in previous papers. 

\bigskip \noindent                    For $\sigma \in \Sio=\{a, b\}^\om$ we denote $\sigma^a$ the $\om$-picture whose first row is the $\om$-word $\sigma$ and whose other rows 
are labelled with the letter $a$. 
 For an \ol~  $L \subseteq \Sio=\{a, b\}^\om$  we  denote  $L^a$  the language of infinite pictures $ \{ \sigma^a \mid \sigma \in L \}$.

\begin{Lem}[\cite{Fin04}] \label{lemTM}
If   $L \subseteq \Sio$ is 
accepted by some Turing machine (in which every run is complete non oscillating) with a B\"uchi acceptance 
condition, then $L^a$ is B\"uchi recognizable by a finite tiling system. 
\end{Lem}

\begin{Lem}[\cite{Fink-tilings}]\label{red1}
There is an injective computable function 
$K$ from $\mathbb{N}$ into $\mathbb{N}$ satisfying the following property. 
\nl  If ${\cal M}_z$ is the non deterministic Turing machine (equipped with a $1'$-acceptance condition) of index $z$, 
and if ${\cal T}_{K(z)}$ is the tiling system 
(equipped with a B\"uchi  acceptance condition) of index $K(z)$, 
then  
$$L({\cal M}_z)^a  = L^B({\cal T}_{K(z)})$$ 
\end{Lem}

\noindent          On the other hand, we can easily see that the cardinality of the $\om$-language $L({\cal M}_z)$ is equal to the cardinality of the 
$\om$-picture language $L({\cal M}_z)^a$.  Thus using the reduction $K$ given in the above lemma we see that: 

$$\{ z \in \mathbb{N} \mid  L({\cal M}_z) \mbox{ is countable } \} \leq_1 \{  z \in \mathbb{N}  \mid  L^B({\cal T}_z) \mbox{ is countable }\}$$
\noindent          Then the completeness part follows from the fact that $\{ z \in \mathbb{N} \mid  L({\cal M}_z) \mbox{ is countable } \}$ is $\Pi_1^1$-complete, 
proved in \cite{FL-TM}. 

\bigskip \noindent                   (2). The proof of Item (2) follows directly from Item (1). 

\bigskip \noindent                   (3). We already know that 
the  set $\{ z \in \mathbb{N} \mid  L^B({\cal T}_z) \mbox{ is infinite} \}$ is  in the class $\Sigma_1^1$. 
Thus  the set $\{ z \in \mathbb{N} \mid  L^B({\cal T}_z) \mbox{ is countably infinite} \}$ is the intersection of a $\Sigma_1^1$-set and of a 
$\Pi_1^1$-set, i.e. it is in the class $D_2(\Sigma_1^1)$.  Using again the reduction $K$ we see that: 
$$\{ z \in \mathbb{N} \mid  L({\cal M}_z) \mbox{ is countably infinite } \} 
\leq_1 \{  z \in \mathbb{N}  \mid  L^B({\cal T}_z) \mbox{ is countably infinite }\}$$
\noindent          It was proved in  \cite{FL-TM} that $\{ z \in \mathbb{N} \mid  L({\cal M}_z) \mbox{ is countably infinite  } \}$ is 
$D_2(\Sigma_1^1)$-complete. Thus the set 
$\{  z \in \mathbb{N}  \mid  L^B({\cal T}_z) \mbox{ is countably infinite }\}$
 is also $D_2(\Sigma_1^1)$-complete. 
\ep 

\bigskip \noindent                   We are now looking at complements of  recognizable languages of infinite pictures. We first state the following result which shows that 
 actually the cardinality of the complement of a 
recognizable language of infinite pictures may depend on the models of set theory. We denote $L^B({\cal T})^-$ the complement 
$\Si^{\om, \om}-L^B({\cal T})$ of a B\"uchi-recognizable language $L^B({\cal T}) \subseteq \Si^{\om, \om}$. 

\begin{The}\label{Complement}
The cardinality of the complement of a B\"uchi-recognizable language of infinite pictures  is not determined by the axiomatic system {\bf ZFC}. 
Indeed there is a B\"uchi tiling system ${\cal T}$ such that: 
\begin{enumerate}
\item There is a model $V_1$ of  {\bf ZFC} in which      $L^B({\cal T})^-$ is countable. 
\item There is a model $V_2$ of  {\bf ZFC} in which     $L^B({\cal T})^-$ has cardinal $2^{\aleph_0}$. 
\item There is a model $V_3$ of  {\bf ZFC} in which      $L^B({\cal T})^-$ has cardinal $\aleph_1$ with 
$\aleph_0<\aleph_1<2^{\aleph_0}$.
\end{enumerate}
\end{The}

\proof  
 Moschovakis gave in   
\cite[page 248]{Moschovakis80} a $\Pi_1^1$-formula $\phi$ defining the  set $\mathcal{C}_1( \Sio)$. 
 Thus its complement 
$\mathcal{C}_1( \Sio)^-=\{a, b\}^\om - \mathcal{C}_1( \Sio)$ is a   $\Si_1^1$-set  defined by the  $\Si_1^1$-formula $\psi=\neg \phi$. 

\bigskip \noindent                    Recall  that one can construct, from the $\Si_1^1$-formula $\psi$ defining $\mathcal{C}_1( \Sio)^-$, 
a B\"uchi Turing machine ${\cal M}$ accepting the $\om$-language 
$\mathcal{C}_1( \Sio)^-$. 

\bigskip \noindent                   On the other hand it is easy to see that the language $\Si^{\om, \om}- (\Sio)^a$ of $\om$-pictures is B\"uchi recognizable. But the class 
$TS(\Si^{\om, \om})$ is closed under finite union, so  we get the following result. 

\begin{Lem}[\cite{Fink-tilings}]\label{lem2'}
If   $L \subseteq \Sio$ is 
accepted by some Turing machine with a B\"uchi acceptance 
condition, then $L^a \cup  [ \Si^{\om, \om} - (\Sio)^a ]$ is B\"uchi recognizable by a finite tiling system. 
\end{Lem}

\noindent          Notice that the constructions are effective and that they can be achieved in an injective way.  Thus we can construct, from the 
B\"uchi Turing machine ${\cal M}$ accepting the $\om$-language 
$\mathcal{C}_1( \Sio)^-$, a B\"uchi  tiling system ${\cal T}$ such that  
$$L^B({\cal T})=L({\cal M})^a \cup  [ \Si^{\om, \om} - (\Sio)^a ].$$
\noindent          It is then easy to see that:
$$L^B({\cal T})^-=(\Sio - L({\cal M}))^a=(\mathcal{C}_1( \Sio))^{a}.$$
\noindent          Thus the cardinality of $L^B({\cal T})^-$ is equal to the cardinality of the $\om$-language $\mathcal{C}_1( \Sio)$, and then 
we can infer the results of the theorem from previous Corollary  \ref{cor1}. 
\ep

\bigskip \noindent                   
We can now use the proof of the above result to prove the following result which shows that natural cardinality problems are  actually 
located at the third level of the analytical hierarchy.

\begin{The}\noindent          
\begin{enumerate}
\ite 
$\{  z \in \mathbb{N}  \mid  L^B({\cal T}_z)^- \mbox{ is finite } \}$ is $\Pi_2^1$-complete. 
\ite 
$\{  z \in \mathbb{N}  \mid  L^B({\cal T}_z)^- \mbox{ is countable } \}$ is in $\Si_3^1 \setminus (\Pi_2^1 \cup \Si_2^1)$. 
\ite
 $\{  z \in \mathbb{N}  \mid  L^B({\cal T}_z)^- \mbox{ is uncountable}  \}$ is in $\Pi_3^1 \setminus (\Pi_2^1 \cup \Si_2^1)$. 
 \end{enumerate}
\end{The}

\proof Item (1) was proved in \cite{Fink-tilings}.  

\bigskip \noindent                   To prove Item (2), we first show that $\{  z \in \mathbb{N}  \mid  L^B({\cal T}_z)^- \mbox{ is countable } \}$ is in the class $\Si_3^1$. 

\bigskip \noindent                   As in \cite{Fink-tilings}, 
 using a recursive bijection $b : (\mathbb{N}-\{0\})^2  \ra \mathbb{N}-\{0\}$,  we can consider an infinite word $\sigma \in \Sio$ 
as a countably infinite family of infinite words over $\Si$ :  the family of $\om$-words 
$(\sigma_i)$ such that for each $i \geq 1$, $\sigma_i$ is defined by $\sigma_i(j)= \sigma(b(i, j))$ for each $j\geq 1$. 
And one can associate with each  $\om$-word $\sigma \in \Sio$ a unique  $\om$-picture $p^\sigma \in \Si^{\om, \om}$ which is simply defined 
by $p^\sigma(i, j)=\sigma(b(i, j))$ for all integers $i, j \geq 1$. 

\bigskip \noindent                   We can now express ``$L^B({\cal T}_z)^- \mbox{ is countable }$" by the formula: 
$$\exists \sigma \in \Sio ~~ \fa p \in \Si^{\om, \om} ~[ ( p \in L^B({\cal T}_z)  ) \mbox{ or } (  \exists i \in \mathbb{N} ~p = p^{\sigma_i} ) ]$$
\noindent          This is a $\Si_3^1$-formula because ``$p \in L^B({\cal T}_z)$", and hence also 
``$[ ( p \in L^B({\cal T}_z)  ) \mbox{ or } (  \exists i \in \mathbb{N} ~p = p^{\sigma_i} ) ]"$,  is expressed by a $\Si_1^1$-formula.

\bigskip \noindent                   We can now prove that $\{  z \in \mathbb{N}  \mid  L^B({\cal T}_z)^- \mbox{ is countable } \}$ is neither in the class $ \Si_2^1$ nor in the 
class $\Pi_2^1$, by using Shoenfield's  Absoluteness  Theorem from Set Theory. 

\bigskip \noindent                   Let   ${\cal T}$ be the  B\"uchi tiling system obtained in  Theorem \ref{Complement} and let $z_0$ be its index 
so that ${\cal T}={\cal T}_{z_0}$.  

\bigskip \noindent                    Assume now  that {\bf V} is  a model of ({\bf ZFC} + $\om_1^{\bf L} < \om_1$).   In the model {\bf V}, by the proofs of  Theorem \ref{Complement}
and  of Corollary  \ref{cor1},   the integer $z_0$ belongs to the set 
$ \{  z \in \mathbb{N}  \mid L^B({\cal T}_z)^- \mbox{ is countable } \}$. 

\bigskip \noindent                   But, by the proofs of  Theorem \ref{Complement}
and  of Corollary  \ref{cor1},  in the inner model     ${\bf L } \subseteq {\bf V }$, the language $L^B({\cal T}_{z_0})^-$ has 
cardinality $2^{\aleph_0}$. 
Thus   the integer $z_0$ does not belong to the set 
$ \{  z \in \mathbb{N}  \mid   L^B({\cal T}_z)^- \mbox{ is countable }  \}$. 

\bigskip \noindent                   On the other hand, Schoenfield's Absoluteness Theorem implies that every  $\Si_2^1$-set (respectively,  $\Pi_2^1$-set)
 is absolute for all inner models of {\rm  (ZFC)}, 
see \cite[page 490]{Jech}.
\nl  In particular, if the set   $ \{  z \in \mathbb{N}  \mid   L^B({\cal T}_z)^- \mbox{ is countable } \}$   was a  $\Si_2^1$-set or a $\Pi_2^1$-set
 then it could not be a different subset of $\mathbb{N}$ in the models  ${\bf  V}$  and   ${\bf L }$ considered above.
Therefore, 
the   set  $ \{  z \in \mathbb{N}  \mid  L^B({\cal T}_z)^- \mbox{ is countable }  \}$   is neither a $\Si_2^1$-set nor a $\Pi_2^1$-set. 

\bigskip \noindent                    Item (3) follows directly from Item (2). 
\ep 

\section{Concluding Remarks}

\noindent           Using the notion of largest effective coanalytic set, 
we have proved in another paper that the topological complexity of a recognizable language of infinite pictures 
is  not determined by the  axiomatic system {\bf ZFC}.  In particular, there is a B\"uchi tiling system $\mathcal{S}$ and models ${\bf V}_1$ and ${\bf V}_2$
of {\bf ZFC} such that: the $\om$-picture language $L(\mathcal{S})$ id Borel in ${\bf V}_1$ but not  in ${\bf V}_2$, \cite{Fin-ICST}. 

\bigskip \noindent                    We have proved in this paper that
$\{  z \in \mathbb{N}  \mid  L^B({\cal T}_z)^- \mbox{ is countable } \}$ is in $\Si_3^1 \setminus (\Pi_2^1 \cup \Si_2^1)$. 
It  remains open  whether this set is actually $\Si_3^1$-complete.

\end{document}